\documentclass[12pt, noinfoline]{imsart}

\usepackage{amsxtra,amssymb,multicol,graphicx,hyperref,color}
\usepackage[normalem]{ulem}
\usepackage{geometry,marginnote}
\usepackage{verbatim,graphicx}

\begin{document} 
\title{WHAT IS ... a Markov basis?}
\author{Sonja Petrovi\'c}
\maketitle
\vskip1in

\begin{figure}[!bh]
  \begin{center}
\includegraphics[scale=.4]{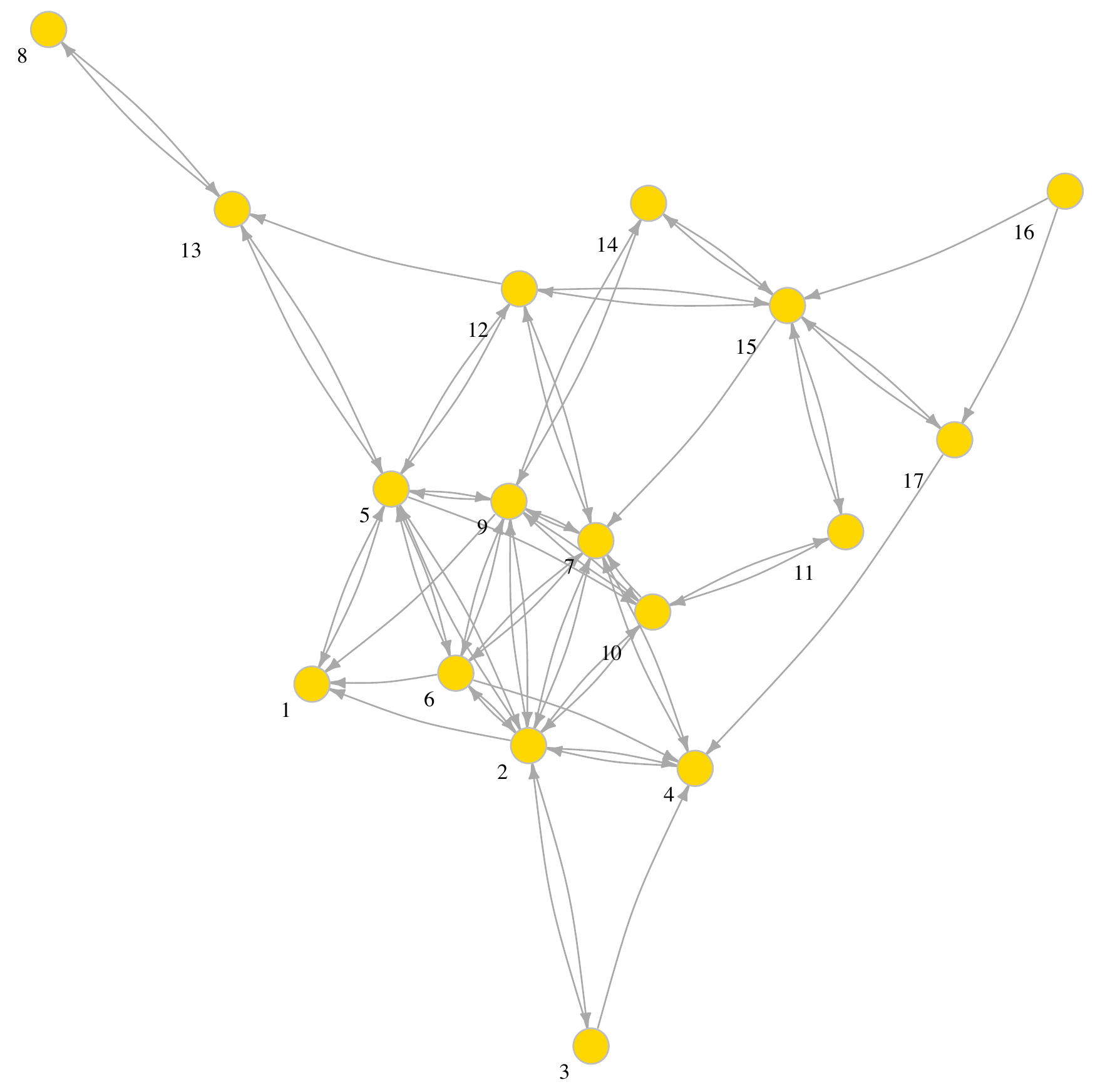}
  \end{center}
 \caption{A graphical summary of the shareholding relationships of those corporations who disclosed them. Nodes in the graph represent the corporations. A directed edge from $i$ to $j$ means company $i$ owns shares in company $j$. }
 \label{fig:nytplot}
\end{figure}

A June 25th, 2014 New York Times article \emph{`Shinzo Abe's Bid to Shake Up Corporate Japan'} by  Hiroko Tabuchi discussed share ownership of 
  $17$  Japanese corporations and cited Prime Minister Abe's claim that they are connected in interlocking ways, owning shares in each other `to create relationships that can protect them from outside interference'. 
  A graphical representation of this relationship is depicted in Figure~\ref{fig:nytplot}. 
Can this claim be verified? How confident are we that it is anything beyond basic intuition? 

Questions like this are at the heart of statistical reasoning: given an observed data set, we wish to find out how surprising it is given some assumption about the world. The fact is, we face such questions on a regular basis: 
 do male faculty members have higher salaries than their female counterparts? We answer them by looking  at the data on salaries and breaking it down by gender, discipline, etc. Our hope that gender is independent of salary translates in a very intuitive way into an expectation: we expect  to see certain distribution of numbers in the salary data; when we do not, we suspect that our assumption of independence could be wrong. We shall 
  put this everyday intuition into a formal framework: the assumption of independence is captured by a statistical model, a family of probability distributions for the salary data that take a specific form;  the observed data is evidence for or against such a model. The evidence carries with it a weight, a confidence level which measures its strength. 
 
 \section*{The algebra behind the intuition}
While the idea of testing a model based on observed data is a simple one in statistics, an exciting development in the 1990's has brought to bear tools from commutative algebra to help solve the problem for a special class of statistical models which we can define using an integer matrix. 
A Markov basis is a set of vectors in the null space of that matrix that allows us to generate synthetic data, starting from the observed, and use the resulting sample to gather evidence against or for the proposed model.

%


Mathematically,  given a statistical model defined by $A\in\mathbb Z^{m\times r}$, 
a \emph{Markov basis for the model} is a set of vectors $\{b_1,\dots,b_n\} \subset \ker_{\mathbb Z} A$ such that  for 
 every pair of vectors $u,v$ for which $Au=Av$,
there exists a choice of basis vectors  satisfying 
\[
	u + b_{i_1} + \ldots + b_{i_N} = v, 
\]
where each partial  sum results in a non-negative vector, that is, 
$
	u+\sum_{j=0}^N b_{i_j}  \geq 0, 
$
componentwise for any $j=1\dots N$. 
Set  ${b_i^{+}}  = \min(0,b_i)$ and ${b_i^{-}}=\max(0,b_i)$, so that each vector in the basis can be written as a difference of nonnegative vectors $b_i=b_i^{+}-b_i^{-}$. 

{\bf Theorem \cite{DS98}.} A  set of vectors is a Markov basis \emph{if and only if} the corresponding set of binomials 
$
	\{ x^{b_i^{+}} - x^{b_i^{+}} \}
$
generates the toric ideal $I_A:=( x^u-x^v : u-v\in\ker_{\mathbb Z}A )$. 

One of the remarkable consequences of this theorem is that the existence of a finite Markov basis for \emph{any} model that can be defined by such a matrix  is  now guaranteed by the Hilbert basis theorem.
Besides this being a lovely mathematical result connecting commutative algebra with statistics, it turns out that Markov bases are  a necessary  tool  for reasoning with  certain types of data such as the Japanese corporate example.

\section*{Formal reasoning with data}
In order to restate our opening questions more formally, 
let us think like statisticians: buying shares is a random event that occurs with some probability. The New York Times  article  suggests that this probability is governed in large part by reciprocity:  if you own shares in my company, I am likely to buy shares in yours. 
If we  can construct a statistical model that produces 
relational data where reciprocity matters, 
the question then becomes whether such a model fits the observed set of relationships.
 In other words, we seek to find out whether such a model can adequately explain how the share ownership data was generated.

The statistics literature provides us with a model that was designed to capture precisely this type of a reciprocal relationship
\cite{HL81, FW81}.   
The model  comes equipped with an integer matrix: it is \emph{log-linear} in form, 
an example of a discrete exponential family, defined in the next section. 
We will test its goodness of fit to the observed data using an \emph{exact conditional test}, defined in the following section,  which requires an understanding of a certain conditional probability distribution.  
Markov bases are a key ingredient to this step. 

Statistical reasoning then proceeds as follows: 
The New York Times article's claim can be restated from the perspective of statistical models as follows: \emph{under the assumption} that the random event of buying shares is governed by reciprocity,  the relationship between Japanese corporations is \emph{not} unusual. 
If this is correct, then the model on dyadic relationships with non-zero reciprocation effect 
  fits the corporate network data, while the one in which reciprocation effect is set to zero does not. 
That is,  the observed data is  extreme/an outlier in the latter case, and not so in the former. 
As there are too many possible share-buying scenarios within 17 companies to which we should compare the observed data, we instead sample from a carefully chosen  reference set, one that makes statistical sense. Markov bases, defined for any log-linear model, are used to design an MCMC sampling scheme of this reference set.

\section*{The $A$ matrix: log-linear models} 

Log-linear models are a class of statistical models for discrete data  for which logarithms of joint probabilities  are captured by a linear map as follows. 
Let $X_1,\dots,X_k$ be discrete random variables with $X_i$ taking values in $[d_i]$. A 
$k$-way \emph{contingency table} $u\in\mathbb Z_{\geq0}^{d_1\times\cdots\times d_k}$ is a nonnegative integer table whose $(i_1,\dots,i_k)$-entry counts the number of times the  event $\{X_1=i_1, \dots, X_k=i_k\}$ occurred; we think of $u$ as a realization of a random table $U$. 
A typical example of the use of such tables is in a cross-classification of items into  $k$ categories (e.g., salary levels by gender). 
Fix an integer matrix $A\in\mathbb Z^{m\times (d_1\cdots d_k)}$ such that $(1,\dots,1)$ is  in its row span\footnote{A normalizing assumption so that all of the  details make sense for probability vectors that must sum to $1$.}. Flatten the data table $u$ into a $d_1\cdots d_k$ column vector. We may interpret the vector $Au$ as a summary of the data table $u$. 

{\bf Definition.} With the above setup, 
the statistical model for which the summary $Au$ suffices to capture the probability of $u$ is called \emph{the log-linear model $\mathcal M_A$ for $k$-way tables associated to the matrix $A$}. It is the family of probability distributions of the following form: 
\[ P_\theta(U=u) = \exp\{\left< Au,\theta\right>-\psi(\theta)\},\]
where $\theta\in\mathbb R^{m}$ is the vector of model parameters and $\psi(\theta)$ is the normalizing constant\footnote{This is simply to ensure the probabilities are nonnegative and sum to $1$.}. The $i$-th entry of the vector $Au$ is called the \emph{minimal sufficient statistic} for the parameter $\theta_i$. The matrix $A$ that computes the sufficient statistics is called the \emph{design matrix} of $\mathcal M_A$. 

Consider the simple example of 
   independence of two discrete random variables, $X$ and $Y$, taking values in $[d_1]$ and $[d_2]$, respectively. 
 Let $\alpha_i=P(X=i)$ and $\beta_i=P(Y=i)$ denote the marginal probabilities of $X$ and $Y$. 
The model $\mathcal M$ of independence postulates that the joint probabilities factor as $ P(X=i,Y=j)=\alpha_i\beta_j$. 
In the language of \cite[Definition 1.1.9]{DSSlectures}, 
 $\mathcal M$ is the toric model associated to $A$, because the monomial joint probabilities $\alpha_i\beta_j$ parameterize a toric variety. 
Data on $X,Y$ can  be organized in a $2$-way table, where the $ij$-entry counts the number of occurrences of the event $\{X=i,Y=j\}$. 
Under the model $\mathcal M$, to know the probability of observing a given data table $u$ it suffices to know the marginal probabilities $\alpha_i$s and $\beta_j$s. 
The corresponding sufficient statistics are marginal counts -- row and column sums -- of the data table $u$.  
As computing these marginals is a  linear operation, it can be presented as a linear map $u\mapsto Au$, where  $A\in\mathbb Z^{(d_1+d_2)\times d_1d_2}$ and $u$ is flattened to a $d_1d_2\times1$ vector. 

 \section*{The weight of the evidence: exact conditional $p$-value}
 What, then, is the  conditional test for whether a log-linear model $\mathcal M_A$ fits the observed data table $u$? 
As we seek to answer whether  $u$ is more-or-less expected under $\mathcal M_A$, 
the test approximates  the \emph{exact conditional $p$-value} of $u$: the probability of a data table being more extreme (less expected) than $u$, conditional on the observed values of the sufficient statistics.
Since sufficient statistics  offer a summary of $u$ that fully capture its probability of occurring under $\mathcal M_A$, 
it is reasonable to condition on the value of $Au$ and explore the  resulting distribution and set of tables. 
The  set $$\mathcal F_A(u) := \{v\in\mathbb Z_{\geq0}^{d_1\times\ldots\times d_k}: Au=Av\}$$ is called \emph{the fiber of $u$ under the model $\mathcal M$}, since it is a fiber of the linear map defined by $A$.  

{\bf Definition.}  A \emph{Markov basis} of the model $\mathcal M_A$ is any set of tables $\mathcal B:= \{b_1,\dots,b_n\}\subset \mathbb Z^{d_1\times\ldots\times d_k}$, called `moves', for which 
\[
	A b_i =0 
\] 
and such that 
for any data table $u\in\mathbb Z_{\geq0}^{d_1\times\ldots\times d_k}$ and 
for any $v\in\mathcal F_A(u)$, there exist $b_{i_1},\dots,b_{i_N}\in \mathcal B$ that can be used to reach $v$ from $u$: 
\[
	u + b_{i_1} + \ldots + b_{i_N} = v
\]
while walking through elements of the fiber:  
\[
	u+\sum_{j=0}^N b_{i_j}  \geq 0 \mbox{, componentwise} 
\]
for any $j=1\dots N$. 

Note that $Au=A(u+b_i)$ means that adding a move $b_i$ to any data table does not change the values of the sufficient statistics, so to remain on the fiber,  we  only need to ensure  that adding a move did not produce negative table entries.

\section*{What's in a basis?}  
The notion of a Markov basis is different (stronger) than that of a basis in linear algebra. Fixing a model and an observed data point results in a fixed conditional distribution of interest.  Think of the finitely many points in this distribution as lying on an integer lattice and  Markov moves as vectors that can be added to a fixed starting point to create a random walk on the lattice.  The set is a basis in the sense that such a random walk is guaranteed to connect all points on the fiber without `stepping outside'.  

Every Markov basis  contains a linear-algebra basis of the null space of $A$.  
  Although the latter can be used to reach all tables in the fiber,  it will generally fail  to satisfy the second condition that each intermediate step is a legal table, since adding one of the basis  elements to some data table may inadvertently make some table entries negative, even while preserving the values of the sufficient statistics. 
To satisfy the non-negativity  condition,  a combination of several null space basis elements may have to be used as a single Markov move in order  to reach or move away from a particular table in the fiber. Algebraically, this can be stated as the fact that a generating set of the toric ideal $I_A$ can be obtained by saturation from the lattice basis ideal defined by $A$; for more details, see \cite[\S 1.3]{DSSlectures}. 

\section*{The algebraic advantage has its challenges} 
Markov bases are one of the two popular ways to sample from the conditional distribution on the fiber (the other is called sequential importance sampling). Fibers $\mathcal F_A(u)$ are generally far too large to enumerate for most reasonably-sized matrices $A$ in practice, and thus exploring them via random walks is a natural alternative. 
The basic idea of Markov bases is therefore quite straightforward, yet it has provided a multitude of open problems over the past two decades.  
The terminology 
  was coined in \cite{DS98} and 
it has become a cornerstone of one area  of algebraic statistics.

Let us revisit the algebraic setup. Looking back to the independence model example, 
let us arrange the joint probabilities of $X$ and $Y$ 
 in a table $p\in[0,1]^{d_1\times d_2}$.  
A probability table $p$ is in the model if and only if it is of rank $1$, or, equivalently, can be written as an outer product of the two marginal probability vectors ($p_{ij}=\alpha_i\beta_j$). Rank-one is of course a determinantal condition: $p_{ij}p_{kl}-p_{il}p_{kj}=0$ for all $i,j,k,l$. This binomial corresponds to a Markov move that replaces the  $1$s in positions $il$ and $kj$ of the table with $1$s in positions $ij$ and $kl$.
 It is one of the defining polynomials of the toric ideal associated to the design matrix of the independence model. Hilbert basis theorem guarantees that  this ideal is finitely generated;  in fact, \emph{any set of generators of this ideal is a Markov basis of the model}. 
The correspondence between bases connecting the fibers and generating sets of toric ideals 
 is often called the Fundamental Theorem of Markov Bases \cite[\S 1.3]{DSSlectures}. 
 Random walks on fibers  constructed using Markov moves 
  come with certain convergence guarantees - if they are used as proposal moves in a Metropolis-Hastings algorithm to sample the fiber, the stationary distribution of that Markov chain will---by design---be precisely the conditional distribution we are interested in sampling.

How complicated are Markov bases? In a wonderful theorem about decomposable graphical models, a special class of models that can be broken into components recursively, Dobra provided a divide-and-conquer strategy to compute Markov bases for all such models on tables of any dimension. Then came the fundamental bad news result of De Loera and Onn: the easiest of non-decomposable models on three-dimensional contingency tables is such that Markov bases can be `as complicated as you can imagine' if two of the dimensions are allowed to grow. (If two dimensions are fixed, then the moves are of bounded complexity. All of these results are summarized in \cite[\S 1.2]{DSSlectures}.)  
So a natural question arises: if the model is not decomposable, how does one compute a Markov basis or verify that a proposed set of moves constitutes one? 
One general strategy is to use the so-called distance-reducing method \cite{AHT}: consider two arbitrary points in the fiber and show that the distance between them can be reduced by applying some of the proposed moves. 

This 
hints at an obvious limitation of Markov bases: \emph{most} of the moves are  unnecessary in many scenarios! Specifically, the bases are \emph{data-independent} by definition: for a fixed $A$,  they connect  the fiber $\mathcal F_A(u)$  for \emph{any} data table $u$, 
 so that for a specific observed data table $u$, many of the moves are inapplicable. 
For example, the following move (right) from the independence model is not applicable to the table on the left:
\[\tiny
	\left[\begin{matrix} 
		2&3&4\\
		0&3&4\\
		0&0&1\\
	\end{matrix} \right]
	+ 
	\left[	\begin{matrix} 
		1&0&-1\\
		0&0&0\\
		-1&0&1\\
	\end{matrix}\right].
\]
Even though it preserves the row and column sums, it produces a negative entry, which is not a `legal' data table that counts the number of occurrences of the event $\{X=3, Y=1\}$. Other restrictions---such as sampling constraints,  maximum on  table entries---make the issue worse. 

The literature offers  a myriad of  results on Markov bases that address these problems: structural results of and complexity bounds for  moves for many classes of models, dynamic algorithms that construct only applicable moves and not an entire basis, larger bases that guarantee connectedness of restricted fibers such as those consisting only of 0/1 tables, etc.

\section*{The $p$-value of shareownership}
Let us finally address Prime Minister Shinzo Abe's aim to diversify the interlocking Japanese corporations. 
The model of interest for this data is also log-linear in form as shown in \cite{FW81}: 
its sufficient statistics 
 are, for each company $i$, the number of companies in which $i$ owns shares, the number of companies that own shares in $i$, and the number of times $i$ reciprocated a shareholding relationship. Computing these counts is a linear operation on the set of companies since it amounts to counting neighbors in the graph above. 
 
We sample the fiber of the observed Japanese corporation data using the dynamic Markov bases implementation for log-linear network models from \cite{GPS}. 
   For each of the $100,000$ sampled data points, we compute a goodness-of-fit statistic---in this case the chi-square statistic---which measures the distance of the data point from what is expected under the model. We do this for two model variants: first when there is a positive reciprocation effect, and second when the reciprocation effect is zero. We take a look at the histogram of this statistic: \emph{the number of times a `more extreme' data point is encountered} is the volume of the histogram to the right of the vertical red line marking the observed value. Since these data are \emph{farther} from the expected value than the observed data, the size of the histogram to the right of the line gives the $p$-value of the data. 
The $p$-values of the data under the models with  non-zero 
   and zero reciprocations are $0.319$ and $0.002$, respectively. 
 
\begin{figure}[!bh]
  \begin{center}
\includegraphics[scale=.35]{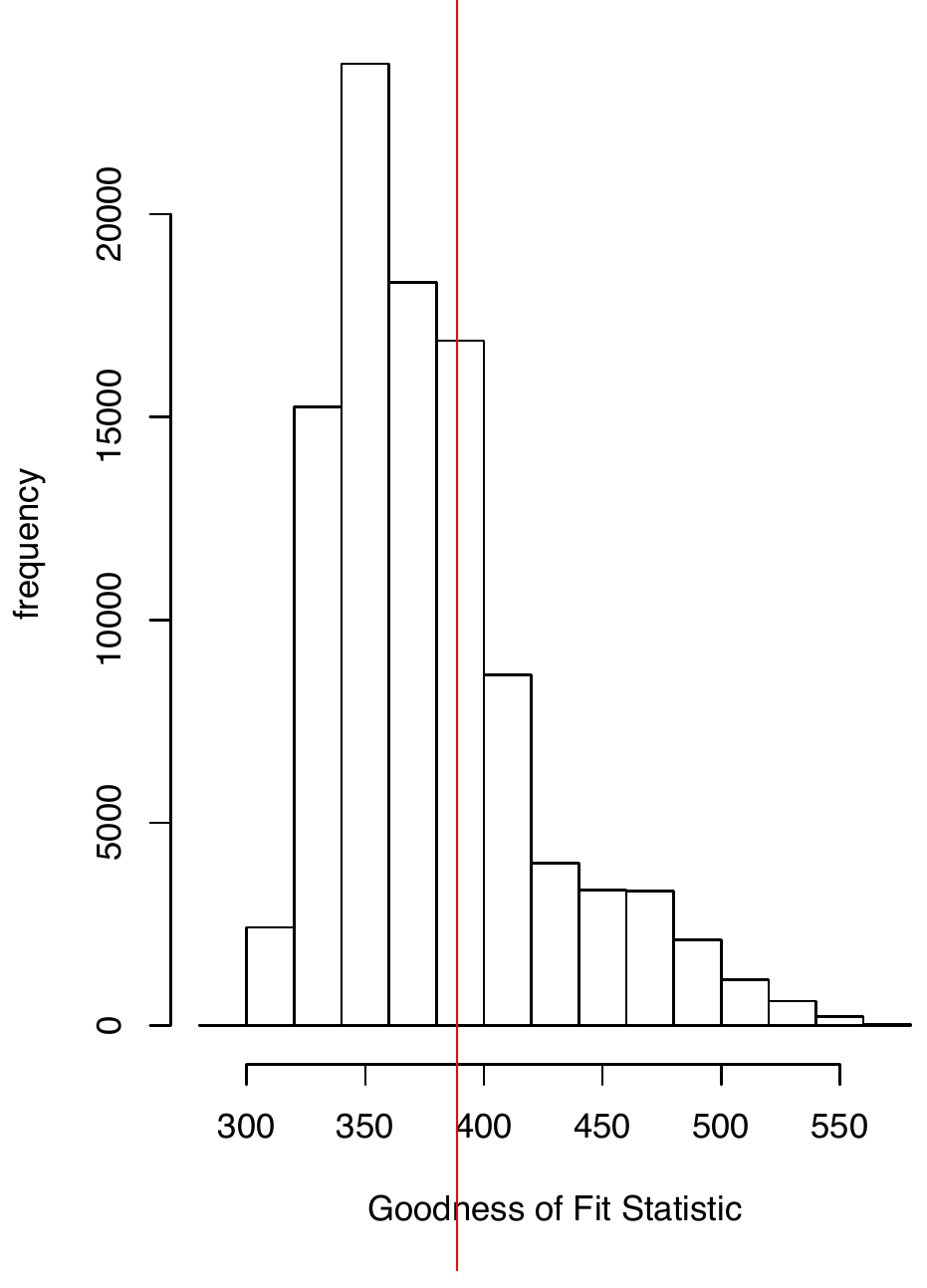}
\includegraphics[scale=.35]{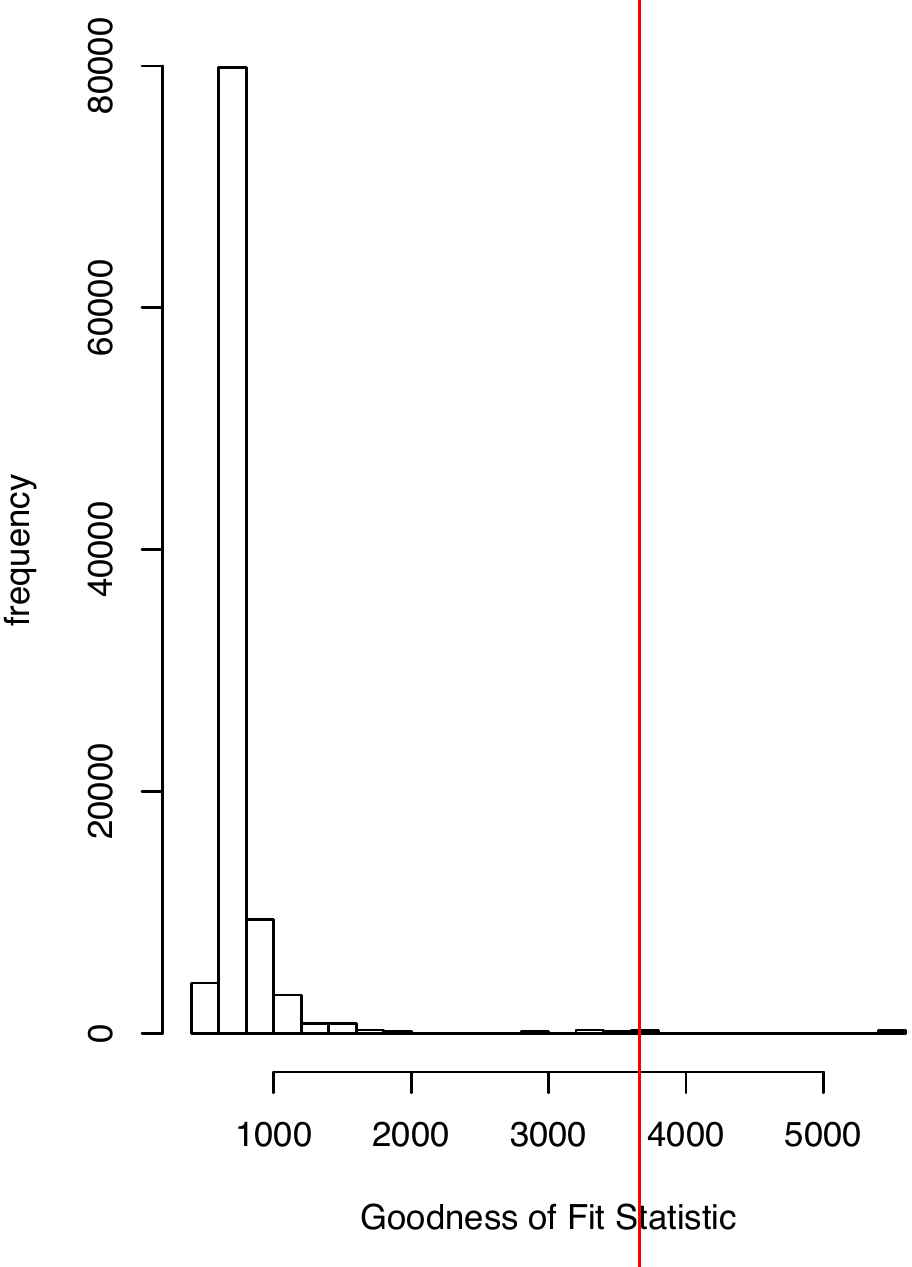}
  \end{center}
 \caption{Histogram of the sampling distribution of the goodness-of-fit statistic. 
Length of each random walk: 100,000 steps. 
Left: model with nonzero reciprocation effect. 
Right: model with no reciprocation effect. 
$p$-values: $p=0.319$ and $p=0.002$, respectively. }
 \label{fig:histograms}
\end{figure}

So, is there statistical evidence to support Prime Minister Shinzo Abe's claim about strong reciprocation effect in interlocking corporate directories? 
Indeed, you may decide that there is, by looking at the histograms in Figure~\ref{fig:histograms}. 



Under the model in which reciprocation effect is present, as many as $31.9\%$ of data points in the sample of $100,000$ are less expected than the observed data, whereas 
under the model with zero reciprocation, that number is merely $0.2\%$. Therefore, the model that sets the reciprocation effect to zero does not fit the data. Perhaps the Prime Minister knew about Markov bases. 

\smallskip
{\bf ACKNOWLEDGEMENT:} {The author is grateful to two PhD students, Lily Silverstein at University of California at Davis and Dane Wilburne at Illinois Institute of Technology, for valuable feedback on a draft of this article in 2018.}


\end{document}